\documentclass[11pt]{article}

\usepackage{amscd,amsmath, amssymb}

\newcommand{\la}{\lambda}

\newcommand{\f}{\varphi}

\newcommand{\CP}[1]{\mathbb{C}P^{#1}}

\numberwithin{equation}{section}

\def\eqref#1{(\ref{#1})}

\newcommand{\arrow}{{\:\longrightarrow\:}}
\newcommand{\Z}{{\Bbb Z}}
\newcommand{\C}{{\Bbb C}}
\newcommand{\R}{{\Bbb R}}

\newcommand{\6}{\partial}
\def\1{\sqrt{-1}\:}
\newcommand{\restrict}[1]{{\left|_{{\phantom{|}\!\!}_{#1}}\right.}}

\newcommand{\calo}{{\cal O}}

\renewcommand{\tilde}{\widetilde}
\renewcommand{\bar}{\overline}
\renewcommand{\phi}{\varphi}
\renewcommand{\epsilon}{\varepsilon}
\renewcommand{\geq}{\geqslant}

\newcounter{Mycounter}[section]
\newcounter{lemma}[section]
\renewcommand{\thelemma}{{Lemma \thesection.\arabic{lemma}}}
\newcommand{\lemma}{%
     \setcounter{lemma}{\value{Mycounter}}
     \refstepcounter{lemma}
     \stepcounter{Mycounter}
     {\noindent \bf \thelemma.\ }}

\newcounter{claim}[section]
\renewcommand{\theclaim}{{Claim \thesection.\arabic{claim}}}
\newcommand{\claim}{%
     \setcounter{claim}{\value{Mycounter}}
     \refstepcounter{claim}
     \stepcounter{Mycounter}
     {\noindent \bf \theclaim.\ }}

\newcounter{sublemma}[section]
\newcounter{corollary}[section]
\newcounter{theorem}[section]
\renewcommand{\thetheorem}{{Theorem \thesection.\arabic{theorem}}}
\newcommand{\theorem}{%
     \setcounter{theorem}{\value{Mycounter}}
     \refstepcounter{theorem}
     \stepcounter{Mycounter}
     {\noindent \bf \thetheorem.\ }}

\newcounter{conjecture}[section]
\newcounter{proposition}[section]
\renewcommand{\theproposition}
       {{Proposition \thesection.\arabic{proposition}}}
\newcommand{\proposition}{%
     \setcounter{proposition}{\value{Mycounter}}
     \refstepcounter{proposition}
     \stepcounter{Mycounter}
     {\noindent \bf \theproposition.\ }}

\newcounter{definition}[section]
\renewcommand{\thedefinition}
       {{Definition~\thesection.\arabic{definition}}}
\newcommand{\definition}{%
     \setcounter{definition}{\value{Mycounter}}
     \refstepcounter{definition}
     \stepcounter{Mycounter}
     {\noindent \bf \thedefinition.\ }}

\newcounter{example}[section]
\newcounter{remark}[section]
\renewcommand{\theremark}{{Remark \thesection.\arabic{remark}}}
\newcommand{\remark}{%
     \setcounter{remark}{\value{Mycounter}}
     \refstepcounter{remark}
     \stepcounter{Mycounter}
     {\noindent \bf \theremark.\ }}

\newcounter{problem}[section]
\newcounter{question}[section]
\makeatletter

\newcommand{\ps@verbit}{%
  \renewcommand{\@oddhead}{%
          \scriptsize
          {Embeddings of Sasakian manifolds}
          \hfil\tiny {L. Ornea and M. Verbitsky, September 15, 2006}}
  \renewcommand{\@evenhead}{\@oddhead}
  \renewcommand{\@oddfoot}{\hfil\thepage\hfil}
  \renewcommand{\@evenfoot}{\@oddfoot}}

\@addtoreset{equation}{section}
\@addtoreset{footnote}{section}
\makeatother

\def\blacksquare{\hbox{\vrule width 5pt height 5pt depth 0pt}}
\def\endproof{\blacksquare}

\begin{document}
\begin{center}
{\LARGE\bf
Embeddings of compact Sasakian manifolds
\\[3mm]
}

Liviu Ornea\footnote{Liviu Ornea is partially supported by  CEx grant no.\ 2-CEx 06-11-22/25.07.2006.} and Misha Verbitsky\footnote{Misha Verbitsky is an EPSRC
advanced
fellow supported by  EPSRC grant
GR/R77773/01.


}

\end{center}

{\small
\hspace{0.15\linewidth}
\begin{minipage}[t]{0.7\linewidth}
{\bf Abstract} \\
Let $M$ be a compact Sasakian manifold. We show that
$M$ admits a CR-embedding into a Sasakian manifold
diffeomorphic to a sphere, and this 
embedding is compatible with the respective Reeb fields. 
We argue that a stronger embedding theorem cannot
be obtained. We use an extension theorem
for K\"ahler geometry: given a compact K\"ahler
manifolds $X\subset Y$, and a K\"ahler form
$\omega_X$ on $X$ which lies in a K\"ahler class
$[\omega]$ of $Y$ restricted to $X$, $\omega_X$ can be
extended to a K\"ahler form $\omega_Y$ on $Y$.
\end{minipage}
}

\tableofcontents


\section{Introduction}
\label{_Intro_Section_}

Sasakian manifolds are in many ways odd-\-di\-men\-sio\-nal
counterparts of K\"ah\-ler manifolds. Quite a few theorems
of K\"ahler geometry have their analogues in Sasakian
geometry.

In K\"ahler geometry one has the Kodaira-Nakano 
embedding theorem stating that a compact 
K\"ahler manifold with integral K\"ahler class is
projective. The embedding given by this theorem is
holomorphic but, in general, non-isometric.

\hfill

Does there exist a similar result in Sasakian geometry?

\hfill

To answer the question, one needs to find a model space
(analogous of the projective space).  One also needs
to find out what is an analogue of holomorphic maps
in Sasakian geometry. 

Sasakian manifolds are equipped with a contact
CR-structure and a CR-holomorphic action of the
corresponding Reeb field. An analogue of a holomorphic
map is obviously a CR-holomorphic immersion
$X\hookrightarrow Y$. One would also require 
that the Reeb field action on $Y$ preserve $X$.

In K\"ahler geometry, the model space for the
Kodaira-Nakano embedding theorem is $\C P^n$. 
However, its analogue
in Sasakian geometry must have continuous
moduli, as we shall now explain. The CR-embeddings 
behave, in this sense, much more like K\"ahler
immersions than the holomorphic immersions.

\hfill

Recall that a Sasakian manifold where the
Reeb field acts with compact orbits is called
{\bf quasi-regular}.

\hfill

Consider a CR-holomorphic immersion of CR-manifolds
$X\hookrightarrow Y$, with $Y$ Sasakian, and $X$ preserved
by the Reeb field. This Reeb action defines a Sasakian
structure on $X$. Indeed, it is positive in the sense
of \cite{bgs}. This is analogous to the fact that
a restriction of a K\"ahler metric to a complex
subvariety is again a K\"ahler metric.
It turns out that any such a map is
compatible with a set of Sasakian structures,
and some of these Sasakian structures are quasi-regular.

\hfill

\claim\label{_q-r-existence_Claim_} 
Let $W_1$, $W_2$ be compact Sasakian manifolds and suppose 
there exists a CR-embedding $W_1\hookrightarrow W_2$ 
which commutes with the flows of the respective Reeb fields 
$\xi_1$, $\xi_2$. Then, there exists a quasi-regular
Sasakian structure on $W_2$ whose Reeb flow preserves $W_1$.

\hfill

\noindent{\bf Proof:} Let $G$ be the group of 
Sasakian automorphisms of $W_2$ which preserve $W_1$.
By assumption, $G$  contains the flow of $\xi_2$
and hence it contains a whole torus in which one can deform 
the initial Reeb flow to a quasi-regular 
one (see \cite{ov2} and also \cite{kami}). \endproof

\hfill

Let $\xi$ be this quasi-regular Reeb field on $W_2$
(its flow preserves $W_1$). It is a positive CR-vector 
field on $W_2$ and its restriction to $W_1$ is still 
CR-holomorphic and positive. 
Hence, according to  \cite{bgs}, $\xi$ is a Reeb field on $W_1$ too. 
Then the orbifolds $W_1/\xi$ and $W_2/\xi$ are 
K\"ahler\footnote{The K\"ahler structure on $W_i/\xi$
is obtained as a restriction of the Levi form.},
and one  has the K\"ahler (thus {\em isometric})
embedding $W_1/\xi\hookrightarrow W_2/\xi$.

\hfill

Now, if such a model space $W_2$ exists, 
the K\"ahler orbifold it fibers on is determined by a choice of a 
quasi-regular Reeb field. Such a Reeb field corresponds to a 
homomorphism $S^1\hookrightarrow T$, where the torus $T$
is the group of isometries
of $W_2$ obtained as the closure of the flow of $\xi_2$
(see \ref{_q-r-existence_Claim_}). But there only exists a
countable number of Lie group homomorphisms of $S^1$ 
into a compact torus. The source
Sasakian manifold, $W_1$ in the above notations,
projects on an arbitrary K\"ahler
orbifold. And one cannot expect to have an isometric
embedding of an arbitrary K\"ahler orbifold into 
a model K\"ahler orbifold, or a countable number
of such orbifolds. Obviously, the model space
must have continuous (and even infinite-dimensional)
moduli.

\hfill

This shows that one cannot aim at a precise analogue of the Kodaira-Nakano
theorem for Sasakian manifolds (and, in particular, Theorem 6.1 in \cite{ov2}
is not correct as stated).\footnote{The argument used in
the proof of \cite[Theorem 6.1]{ov2} is essentially correct, but
only a contact immersion is constructed.}

Here we present the best possible result in this direction:

\hfill

\theorem \label{_main_embe_Theorem_}
Any compact
Sasakian manifold admits a CR-embedding into a Sasakian
manifold $S$ which is diffeomorphic to a sphere, and
such that the corresponding K\"ahler manifold 
$S\times \R^{> 0}$ is biholomorphic to
$\C^n\backslash 0$.

\hfill

We prove \ref{_main_embe_Theorem_} in Section \ref{_embe_Section_}.

\section{Background on Sasakian and Vaisman manifolds}
\label{_Sasa_Section_}

We present the minimum necessary facts about Sasakian and
Vaisman manifolds. For details and examples, see
\cite{_BG:Book_, bgs, drag, ov1, ov2, ov3, ov4, _Vaisman_}.

\hfill

\definition A {\bf Sasakian manifold} is a Riemannian
manifold $(W, g)$ such that the cone metric $t^2g+dt^2$ on
$W\times \R^{>0}$ is K\"ahler with respect to a
dilatation-invariant complex structure.

\hfill

The conic K\"ahler metric has a global potential $\f_W=t^2$.
Clearly, $W$, identified with the slice at $t=1$, is the
level set $\f_W^{-1}(1)$. As such, $W$ has a subjacent
strictly-pseudoconvex structure, the Reeb field being
$\displaystyle\xi=J\frac{d}{dt}$ (see \emph{e.g.}
\cite{ov4}). It can be shown that the Reeb field is
Killing, unitary and all the plane sections passing
through it have curvature $1$:

\begin{equation}\label{xi}
R(X,Y)\xi=g(\xi, Y)X-g(X,Y)\xi.
\end{equation}

A pair $(g,\xi)$, with $\xi$ Killing, unit and
satisfying \eqref{xi} is called a {\bf Sasakian structure}
on $W$.

The foliation generated by $\xi$ has a transverse K\"ahler
structure. If the orbits of $\xi$ are closed, the Sasakian
structure is called {\bf quasi-regular}. In this case, the
Reeb field generates a locally free $S^1$-action such that
the leaf space is an orbifold and the transverse K\"ahler
structure projects to it. Hence, a quasi-regular Sasakian
manifold is a principal $S^1$-bundle over a K\"ahler
orbifold.

\hfill

\remark \label{example} 
All quasi-regular
Sasakian manifolds can be obtained as follows
(see \cite[Example 1.9]{ov4}). Let $X$ be a projective
orbifold with only quotient singularities and let $L$ be
an ample Hermitian line bundle over $X$ with positive
curvature. Then $\mathrm{Tot}(L^*\setminus\{0\})$, the
space of non-zero vectors in the dual of $L$, has a
K\"ahler structure given by the potential
$\f(v)=\mid v\mid^2$. If $W:=\f^{-1}(1)$ is smooth, it 
is a Sasakian manifold which is a principal $S^1$-bundle 
over $X$.

\hfill

While there are examples of Sasakian structures which are not quasi-regular (even Sasaki-Einstein ones, see \emph{e.g.}
\cite{_GMSW1_}), it is known that a Sasakian structure on a compact manifold can always be deformed to a quasi-regular one (see \cite{ov2}). The trick is to deform the Reeb field in the torus generated by its flow.

\hfill

Sasakian geometry is closely related to Vaisman geometry and our proof exploits this link. We confine to the compact case for which the following definition works (see \cite{ov1}):

\hfill

\definition A {\bf Vaisman manifold} is a Hermitian manifold $(M,J,h)$ admitting a cover which is a cone $W\times \R^{>0}$ endowed with a conic K\"ahler metric globally conformal with the lift of $h$ and with respect to which the deck group acts by holomorphic homotheties.

\hfill

The K\"ahler metric in the definition has a global potential given by $\psi(m,t)=t^2$ (see \cite{_Verbitsky:LCHK_}).

The Reeb field of $W$ and its complex conjugate 
together give rise to a 1-dimensional
holomorphic foliation $\mathcal{F}$ on $M$ which defines a
transversal K\"ahler structure on a Vaisman manifold. 
When this foliation has compact leaves, the Vaisman structure is called
quasi-regular and the leaf-space $M/\mathcal{F}$ is a
K\"ahler orbifold.

It is known that the product of a compact Sasakian manifold with $S^1$ has a natural Vaisman structure.

A typical example of Vaisman manifold is the diagonal Hopf manifold $H^n(\la_1,\ldots,\la_n):=\C^n\setminus\{0\}/\Z\cong S^{2n-1}\times S^1$, where $\Z$ is represented by the group
generated by $(z_i)\mapsto (\la_iz_i)$, for some (ordered) nonzero
complex numbers $\lambda_i$ of module $<1$. Hence, the Vaisman structure is here associated to a deformation of the standard Sasakian structure of the sphere. The orbifold $H^n(\la_1,\ldots,\la_n)/\mathcal{F}$ is in this case the weighted projective space ${\CP{n}}{(\la_1,\ldots,\la_n)}$.

 The Hopf manifold plays in Vaisman geometry the role of the projective space from K\"ahler geometry:

\hfill

\theorem \cite{ov3}\label{embedding_vaisman} 
Every compact Vaisman manifold admits a holomorphic 
embedding in a diagonal Hopf manifold.

\hfill

Now, starting  with a compact quasi-regular Sasakian
 manifold, we 
consider a compact quasi-regular Vaisman manifold
$W\times S^1$. Clearly, leaf-spaces $W/\xi$ and
$M/\mathcal{F}$ are naturally identified. In particular, the weighted
projective space can be viewed as $H^n(\la_1,\ldots,\la_n)/\mathcal{F}$ or as
basis of the Hopf fibration $S^{2n+1}/\xi$. Moreover,  it
can be shown (\cite{ov2, ov3}) that the above holomorphic
embedding descends to a holomorphic embedding of the
orbifold $W/\xi$ into the weighted projective space
${\CP{n}}{(\la_1,\ldots,\la_n)}$. 

\hfill

On the other hand, any Vaisman structure has a subjacent
locally conformal K\"ahler structure, \emph{i.e.} there
exists an open cover such that the restriction of the
Vaisman metric to each open set of the cover is conformal
to some K\"ahler metric. Hence, as in conformal geometry,
one is led to consider the {\bf weight bundle}
$L'\rightarrow M$ which is the trivial line-bundle
associated to the frame bundle by the representation
$\mathrm{GL}(2n,\R)\mapsto \mid
\det(A)\mid^{\frac{1}{2n}}$, where $n=\dim_\C M$. We proved in \cite{ov2} that
for any compact quasi-regular Vaisman manifold, the
push-forward to $M/\mathcal{F}$ of the complexified of
$L'$ is ample. Using \ref{example}, this gives a way of
constructing quasi-regular compact Sasakian manifolds out
of quasi-regular compact Vaisman ones.


\section{The embedding theorem}
\label{_embe_Section_}


Let $W$ be a compact Sasakian manifold, $M=W\times S^1$ 
the associated Vaisman manifold, and 
$\Psi:M\hookrightarrow H^n(\la_1,\ldots,\la_n)$ the
holomorphic embedding into the diagonal Hopf manifold
described in \ref{embedding_vaisman}.

According to \cite[Example 1.9]{ov4}, the
K\"ahler cone $\tilde M=W\times \R^{>0}$ of $M$
is equipped with a K\"ahler potential $\phi_W$ such that 
$W:=W\times \{1\}=\phi_W^{-1}(1)$. Then $\Psi$ naturally lifts 
to a holomorphic embedding
$\tilde M\hookrightarrow \widetilde {H^n}(\la_1,\ldots,\la_n)=\C^n\setminus \{0\}$.

Now assume that the potential $\phi_W$ is the restriction 
of a potential $\phi_H$ on $\C^n\setminus \{0\}$. Then we have a CR
embedding $W\hookrightarrow \phi_H^{-1}(\lambda)$.
The following lemma implies that this space
is diffeomorphic to a sphere.

\hfill

\lemma\label{_level_set_diffeo_Lemma_}
Let $\phi$ be a positive K\"ahler potential on $V= \C^n\setminus\{0\}$,
satisfying $X_A(\phi) = 2\phi$ for a vector
field $X_A$ on $V$, $X_A\restrict v = A(v)$, where $A$
is a linear operator on $V$. Then the level sets
$\phi^{-1}(\lambda)$ are diffeomorphic to a sphere,
for all positive $\lambda\in \R$.

\hfill

\noindent{\bf Proof.}
Since $X_A(\phi) = 2\phi$, $\exp(A)$ acts on $V$ as a homothety,
with respect to the K\"ahler metric defined by $\1\6\bar\6\phi$.
Therefore, $\exp(-A)$ is a contraction. 

Consider the 1-dimensional foliation ${\cal X}_A$ on
$V$ generated by the vector field $X_A$.
Let $r(t,x):= \exp(tA)x$, $x\in V$,
$t\in \R$. Denote by $S_x$ the image
of the map $r(\cdot, x):\; \R \arrow V$.
Clearly, the sets $S_x$ are leaves of ${\cal X}_A$.
The leaf space $V/{\cal X}_A$ is
naturally isomorphic to $S^{2n-1}$ (this isomorphism is apparent,
because $\exp(-A)$ is a contraction).

On $S_x$,
the function $\phi$ becomes exponential after a coordinate
change. Therefore, $S_x$ intersects with
$\phi^{-1}(\lambda)$  precisely once.
This defines a diffeomorphism
\[
  \phi^{-1}(\lambda)\arrow V/ {\cal X}_A= S^{2n-1}.
\]
\endproof

\hfill

We thus have a CR embedding
$W\hookrightarrow \phi_H^{-1}(\lambda)$ in a Sasakian
manifold diffeomorphic to a sphere. 

\hfill

\remark
This diffeomorphism is
not a CR-equivalence. Indeed, CR-isomorphic pseudoconvex
hypersurfaces  in $\C^n$ bound domains which are
biholomorphically equivalent (see e.g. \cite{ov4}). 
On the other hand, the sphere $\phi_H^{-1}(\lambda)$ is built 
from a homogeneous potential which can be chosen arbitrarily.
In \cite{ov4} it is shown that starting from 
an arbitrary strictly pseudoconvex, $S^1$-symmetric
hypersurface $V\subset \C^n$, we can find a K\"ahler
potential for which $V$ is a level set. Of course, a domain 
which is bound by such a hypersurface is not necessarily
biholomorphic to a Euclidean ball.

\hfill

From \ref{_level_set_diffeo_Lemma_}, 
we obtain that the embedding problem for Sasakian manifolds is reduced to an
extension problem for plurisubharmonic functions which we now prove.

\hfill

\theorem Let $\Psi:W\times S^1\hookrightarrow {H^n}(\la_1,\ldots,\la_n)$ be the
holomorphic embedding of Vaisman manifolds described above.
Then there exists a K\"ahler potential on $\C^n\setminus\{0\}$
extending the potential $\phi:W\times \R^{>0}\rightarrow \R$ 
and the restriction of  $\Psi$ to the respective level
sets is a CR-embedding.

\hfill

\noindent{\bf Proof.} Replace the Reeb flow on $W$ by 
a quasi-regular one as in
\ref{_q-r-existence_Claim_}.
Consider the holomorphic embedding 
 $\Psi:W\times S^1\hookrightarrow {H^n}(\la_1,\ldots,\la_n)$
of quasi-regular Vaisman manifolds constructed as
in \cite{ov3}.
Let $\tilde M\rightarrow M$ be a  cone covering
of $M= W\times S^1$. We assume that $\tilde M$ is
endowed with a global K\"ahler potential given by the
square length of the Lee field. As in \cite{ov4}, 
let $L\rightarrow X$ be an ample $\C^*$ Hermitian bundle  
over a projective orbifold $X$ such that 
$\tilde M =\mathrm{Tot}(L^*\setminus\{0\})$.

{}From the definitions it is apparent
that the homogeneous of weight $2$ K\"ahler potentials on
$\tilde M$ are in one-to-one correspondence
with the positive Hermitian structures on $L$.

As $X$ is projective, let $f:X\rightarrow \C P^n(\la_1,\ldots,\la_n)$ be a 
holomorphic embedding (in the category of orbifolds, see \cite{baily}).
Of course, $f$ is the projection of  the holomorphic embedding $\Psi$.

The bundle $L$ is the push-forward of the complexification of the 
weight bundle of the Vaisman manifold $W\times S^1$ 
(see \cite{ov2}). As such, one sees that $L$ is the pull-back of the 
tautological bundle of the weighted projective 
space: $L=f^*O(\la_1,\ldots,\la_n)$. 
Hence, to finish the proof, it will be enough to 
show that the Hermitian metric of $L\rightarrow X$ 
is the restriction of a Hermitian metric 
in $O(\la_1,\ldots,\la_n)\rightarrow  \C P^n(\la_1,\ldots,\la_n)$. 
But this will follow if we show that the metric 
(equivalently, the K\"ahler form) of $X$ is the restriction of a 
metric on $\C P^n(\la_1,\ldots,\la_n)$
preserving the given cohomology class,
which is proven in Section \ref{_Kahler_Section_}.


\section{K\"ahler metrics extended from a submanifold}
\label{_Kahler_Section_}


\subsection{K\"ahler metrics and plurisubharmonic functions}

The sketch of the proof of the following theorem was communicated
to us by J.-P. Demailly in an e-mail.

\hfill

\theorem\label{_Kah_exte_Theorem_}
Let $(M, \omega)$ be a compact K\"ahler manifold,
and $Z\subset M$ a closed complex submanifold.
Denote by $[\omega]\in H^2(M)$ the K\"ahler class of $M$.
Consider a K\"ahler form $\omega_0$ on $Z$ such that
its K\"ahler class coincides with the restriction
$[\omega]\restrict Z$. Then there exists a
K\"ahler form $\omega'$ on $M$ in the same K\"ahler
class as $\omega$, such that $\omega'\restrict
Z=\omega_0$.

\hfill

\noindent{\bf Proof.} The proof is based on two theorems,
which are found in \cite{_Demailly_1982_}.
The first one is existence of regularized maximum 
of two plurisubharmonic functions.

\hfill

\proposition\label{_regulari_max_Proposition_}
Let $\phi$ and $\psi$ be two smooth plurisubharmonic
functions on a complex manifold $M$. Then for any
$\delta>0$ there exists a smooth  plurisubharmonic
function $\max_\delta(\phi, \psi)$ which is equal to
$\max(\phi, \psi)$ unless $|\phi-\psi| <\delta$.
Moreover, $\max_\delta(\phi, \psi)$ is strictly
plurisubharmonic if $\phi$ and $\psi$ are strictly
plurisubharmonic.

\hfill

\noindent{\bf Proof.} Let $\lambda:\; \R^2 \arrow \R$ be a convex
function which is non-decreasing in each variable.
Then $\lambda(\phi, \psi)$ is also plurisubharmonic.
Taking a smooth convex approximation $\lambda_\delta$
of $\lambda(x, y) = \max(x,y)$, we obtain 
$\max_\delta(\phi, \psi):= \lambda_\delta(\phi, \psi)$.
\endproof

\hfill

The second one is existence of a plurisubharmonic function
with a logarithmic pole in any given complex variety.

\hfill

\proposition\label{_fu_with_loga_poles_Proposition_}
Let $(M, \omega)$ be a compact K\"ahler manifold
and $Z\subset M$ a complex analytic subvariety.
Then there exists a  function $\phi$
on $M$ which is smooth outside $Z$, satisfies
\begin{equation} \label{_>-C_omega_Equation_}
\1\6\bar\6\phi > - C\omega,
\end{equation} 
for some positive constant
$C$, and has logarithmic
poles along $Z$; that is, locally along $Z$ we have
\[
\phi(z)\sim \log \sum |g_k(z)| + O(1),
\]
where $g_k$ is a local system of generators
of the ideal sheaf of $Z$ in $\calo_M$.

\hfill

\noindent{\bf Proof.} See Lemma 2.1 of \cite{_Demailly_Paun_}.
\endproof

\hfill

Using these two results, we can easily prove
\ref{_Kah_exte_Theorem_}. Write the K\"ahler form
$\omega_0$ as 
\[
\omega_0 = \omega\restrict Z + \1 \6\bar\6 u,
\]
where $u$ is a smooth function on $Z$. Extending
$u$ to $M$ arbitrarily and adding the distance term
$C_1d(z,Z)^2$, we may assume that $u$ is defined 
over $M$ and that $(1-\epsilon)\omega+i\6\bar\6 u$ is strictly 
positive near $Z$, for some positive $\epsilon<\frac 1 2$.
Consider a function $\phi$ with logarithmic poles
in $Z$ constructed in \ref{_fu_with_loga_poles_Proposition_}.
Near $Z$ the function $\phi$ is very small,
hence $\psi:=\max_\delta(\frac \epsilon {C} \phi+A, u)$ is equal to $u$
in some open neighbourhood of $Z$, for any constant $A\in \R$.
Therefore, $\omega':=\omega + \1 \6\bar\6 \psi$ satisfies
$\omega'\restrict Z=\omega_0$. 
Let $U$ be a neighbourhood of $Z$ where
$\1\6\bar\6 u+(1-\epsilon)\omega$ is positive.
The same argument
as used to prove \ref{_regulari_max_Proposition_}
also gives an estimate of $\6\bar\6\psi\restrict U$:
\begin{equation}\label{_6+bar+_6psi_via_u_Equation_}
\1 \6\bar\6\psi\restrict U \geq - (1 - \epsilon)\omega
\end{equation}
(we use $\1\frac \epsilon {C}\6\bar\6\phi\geq -\epsilon \omega$, which
comes from \eqref{_>-C_omega_Equation_}, and 
$\1 \6\bar\6 u\restrict U \geq - (1 - \epsilon)\omega$,
which is essentially a definition of $U$; also, $\epsilon
<1-\epsilon$, because we chose $\epsilon <\frac 1 2$).
Choosing $A$ sufficiently big, we may assume
that $\psi=\phi+A$ outside of $U$. 
Then, outside of $U$, $\omega'$ is positive,
because
\[ \omega'\restrict{M\backslash U}
   =\1 \6\bar\6 \psi+\omega=
   \1\frac \epsilon {C} \6\bar\6 \phi+\omega > (1-\epsilon)\omega
\] 
which is positive.
In $U$ the form $\omega'$ is positive, as follows from
\eqref{_6+bar+_6psi_via_u_Equation_}:
\[
\omega+ \1 \6\bar\6 \psi\geq\omega - (1 - \epsilon)\omega =\epsilon\omega.
\]
We proved \ref{_Kah_exte_Theorem_}.
\endproof

\subsection{K\"ahler forms on orbifolds}

In order to prove \ref{_main_embe_Theorem_},
as indicated in Section \ref{_embe_Section_},
we need to apply the K\"ahler extension result
(\ref{_Kah_exte_Theorem_}) to  weighted 
projective spaces, which are  orbifolds. 
The proof of \ref{_Kah_exte_Theorem_} is valid 
in the orbifold situation, as we shall 
explain now.

{\bf A complex orbifold} is a topological space
covered by open sets of form $B/G$, where $B$ is 
an open ball in $\C^n$, and $G$ a finite group
holomorphically acting on $B$ (different, generally
speaking, for different elements of the covering), 
with the  gluing map holomorphic and equivariant.
It is {\bf K\"ahler} if these open balls
are equipped with a compatible K\"ahler metric, and 
the groups $G$ acts by isometries.

Results about differential forms and the Hodge theory 
 are carried over to orbifolds 
word-by-word. In particular, the $\6\bar\6$-lemma
is valid: given an exact $(p,q)$-form 
$\eta$ on a compact K\"ahler orbifold,
we may always find $\eta_1$ such that
$\eta=\6\bar\6\eta_1$.

The proof of \ref{_Kah_exte_Theorem_}
for orbifolds is done in the same way. We start 
with a smooth (in the sense of orbifolds)
embedding $Z\subset M$, a K\"ahler form
$\omega_0$ on $Z$, and a K\"ahler form
$\omega$ on $M$ with $\omega\restrict Z- \omega_0$
exact. Using $\6\bar\6$-lemma, we write
\[
\omega_0 = \omega\restrict Z + \1 \6\bar\6 u.
\]
The rest of the argument is local, and is repeated
verbatim.

\hfill

\noindent {\bf Acknowledgements:} We are grateful to
J.-P. Demailly for invaluable help, Charles P. Boyer 
for insightful comments and to the referee for 
useful remarks.

\hfill

{\small

\noindent {\sc Liviu Ornea\\
University of Bucharest, Faculty of Mathematics, \\14
Academiei str., 70109 Bucharest, Romania.}\\
\tt Liviu.Ornea@imar.ro, \ \ lornea@gta.math.unibuc.ro

\hfill

\noindent {\sc Misha Verbitsky\\
University of Glasgow, Department of Mathematics, \\15
  University Gardens, Glasgow, Scotland.}\\
{\sc  Institute of Theoretical and
Experimental Physics \\
B. Cheremushkinskaya, 25, Moscow, 117259, Russia }\\
\tt verbit@maths.gla.ac.uk, \ \  verbit@mccme.ru
}

\end{document}